\documentclass[11pt]{article}
\usepackage[margin=1in]{geometry}
\usepackage{amsmath,amssymb,amsthm,mathtools}
\usepackage{booktabs,array,enumitem}
\usepackage{xcolor,graphicx,tikz}
\usetikzlibrary{arrows.meta,positioning,calc}
\usepackage[hidelinks]{hyperref}

\theoremstyle{plain}
\newtheorem{theorem}{Theorem}
\newtheorem{lemma}[theorem]{Lemma}
\newtheorem{corollary}[theorem]{Corollary}
\newtheorem{proposition}[theorem]{Proposition}
\theoremstyle{definition}
\newtheorem{definition}[theorem]{Definition}
\newtheorem{example}[theorem]{Example}
\newtheorem{problem}[theorem]{Problem}

\theoremstyle{remark}
\newtheorem{remark}[theorem]{Remark}
\newcommand{\PC}{\mathcal P}
\newcommand{\dist}{\triangle}
\tikzset{vertex/.style={circle,draw,minimum size=6mm,inner sep=1pt},code/.style={vertex,fill=blue!65!black,text=white},edge/.style={thick},failed/.style={thick,red,dashed}}
\title{Maximum Edge Repair Costs of Perfect Codes}
\author{Leonard M. Paleta\\
\small Department of Mathematics and Statistics\\
\small University of Southern Mindanao\\
\small \href{mailto:lmpaleta@usm.edu.ph}{\texttt{lmpaleta@usm.edu.ph}}}
\date{}

\begin{document}
\maketitle
\begin{abstract}
We study the worst-case number of vertex-membership changes needed to repair a perfect code after one edge deletion. We give general bounds and parity constraints, characterize repairs of costs one and two, determine exact maximum costs for paths, cycles, and complete graphs, and characterize all repaired perfect codes in regular graphs of degree at least two. For vertex-transitive graphs of such degree, finiteness follows from bondage bounds of Huang and Xu. We give a direct proof that, for every failed edge, an automorphic image of any specified initial code survives.
\end{abstract}
\section{Introduction}

Suppose that active sensors occupy selected vertices of a network and that each location must be covered by exactly one sensor, either at that location or at a neighboring vertex. When a link fails, restoring exact coverage may require activating some sensors and deactivating others. How many such changes can one link failure force, and how much does the answer depend on the initial placement?

Let $\Gamma$ be a finite simple undirected graph with vertex set $V(\Gamma)$, edge set $E(\Gamma)$, order $n=|V(\Gamma)|$, and size $m=|E(\Gamma)|$. We follow Diestel~\cite{Diestel2017} for graph terminology and Haynes, Hedetniemi, and Slater~\cite{HHS1998} for domination terminology. The open neighborhood $N_\Gamma(v)$ consists of the vertices adjacent to $v$, and the closed neighborhood is $N_\Gamma[v]=N_\Gamma(v)\cup\{v\}$. The degree of $v$ is $\deg_\Gamma(v)=|N_\Gamma(v)|$. A vertex of degree zero is isolated, and a vertex of degree one is a leaf. An edge incident with a leaf is a pendant edge. We write $P_n$, $C_n$, and $K_n$ for the path, cycle, and complete graph on $n$ vertices, respectively.

A set $S\subseteq V(\Gamma)$ is a \emph{dominating set} if $N_\Gamma[v]\cap S\ne\varnothing$ for every vertex $v$. The \emph{domination number} $\gamma(\Gamma)$ is the minimum cardinality of a dominating set. A \emph{perfect code}, or \emph{efficient dominating set}, is a set $C$ satisfying $|N_\Gamma[v]\cap C|=1$ for every $v$. Equivalently, the closed neighborhoods indexed by the vertices of $C$ partition $V(\Gamma)$. In particular, $C$ is independent, meaning that no two of its vertices are adjacent. Write $\PC(\Gamma)$ for the family of perfect codes of $\Gamma$. Biggs~\cite{Biggs1973} developed the coding-theoretic viewpoint for distance-transitive graphs. The structure and existence of efficient dominating sets have also been studied in the domination literature~\cite{LivingstonStout1990,BangeEtAl1996,HaynesHenningHandbook}.

Edge failures have traditionally been studied through their effect on the domination number. The \emph{bondage number} $b(\Gamma)$ is the minimum number of edges whose deletion increases $\gamma(\Gamma)$~\cite{FinkEtAl1990}. Related work examines domination under edge deletion, vertex deletion, and edge addition~\cite{HaynesHenning2003,XuSurvey}. A graph is $d$-regular if every vertex has degree $d$. An automorphism is a permutation of the vertices preserving adjacency in both directions, and a graph is vertex-transitive if an automorphism maps $u$ to $v$ for every pair of vertices $u,v$. Huang and Xu~\cite{HuangXu2008} obtained bondage bounds for regular and vertex-transitive graphs. Wang and Zhang~\cite{WangZhang2023} studied perfect codes in vertex-transitive graphs. Even when one deletion preserves the domination number, the original placement can require extensive changes. On $C_{3q}$, for example, deleting an edge incident with a code vertex preserves the domination number $q$ but forces $2q$ membership changes (Theorem~\ref{t:cycle}).

To measure these changes, assume throughout that $E(\Gamma)\ne\varnothing$ and $\PC(\Gamma)\ne\varnothing$. For an edge $e$, let $\Gamma-e$ denote the graph obtained by deleting $e$ and retaining all vertices. Given an initial code $C\in\PC(\Gamma)$, a \emph{repair} after deletion of $e$ is any $D\in\PC(\Gamma-e)$. The symmetric difference $C\dist D=(C\setminus D)\cup(D\setminus C)$ records the vertices whose membership changes. The \emph{single-edge repair cost} is
\[
r(\Gamma,C;e)=\min_{D\in\PC(\Gamma-e)}|C\dist D|,\qquad \min\varnothing=\infty.
\]
A repair attaining this minimum is optimal. Each activation and each deactivation costs one unit, so replacing one sensor by another costs two. Exact coverage is required of the final placement.

The \emph{maximum edge repair cost} of $C$ is its worst-case cost over single-edge failures,
\[
R(\Gamma,C)=\max_{e\in E(\Gamma)}r(\Gamma,C;e).
\]
The initial placement can affect this cost even though every perfect code has cardinality $\gamma(\Gamma)$ (Proposition~\ref{p:gamma}). Accordingly, define
\[
R_{\min}(\Gamma)=\min_{C\in\PC(\Gamma)}R(\Gamma,C),
\qquad
R_{\max}(\Gamma)=\max_{C\in\PC(\Gamma)}R(\Gamma,C).
\]
The first quantity gives the smallest worst-case budget achievable by choosing the initial placement. The second gives the budget sufficient for every initial placement. A graph is \emph{fully edge-repairable} if $\PC(\Gamma-e)\ne\varnothing$ for every edge $e$. In this case both budgets are finite (Proposition~\ref{p:finite}), and their difference, the \emph{repair-cost spread} $\sigma_R(\Gamma)=R_{\max}(\Gamma)-R_{\min}(\Gamma)$, measures the largest reduction in worst-case cost achievable through the choice of initial code.

\subsection{Elementary repairs}

The following lemma identifies exactly when a deletion requires at most one membership change.
\begin{lemma}\label{l:elementary-repair}
Let $C\in\PC(\Gamma)$ and $e=uv\in E(\Gamma)$. Then
\[
r(\Gamma,C;e)=0\quad\Longleftrightarrow\quad \{u,v\}\cap C=\varnothing.
\]
If $u\in C$, then $v\notin C$ and
\[
r(\Gamma,C;uv)=1\quad\Longleftrightarrow\quad \deg_\Gamma(v)=1.
\]
\end{lemma}
\begin{proof}
Put $H=\Gamma-e$. Cost zero means that $C$ itself is a perfect code of $H$. If neither endpoint belongs to $C$, deleting $e$ leaves every closed-neighborhood intersection with $C$ unchanged. Otherwise, relabel the endpoints so that $u\in C$. Independence gives $v\notin C$. In $H$, the vertex $v$ loses its unique code neighbor, while every other vertex remains covered exactly once. This proves the first assertion.

A repair of cost one must remove or add a single vertex. Removing a vertex leaves $v$ uncovered. Adding any $x\ne v$ outside $C$ covers $x$ twice, since $x$ was already covered. The only candidate is therefore $C\cup\{v\}$. It covers $v$ exactly once, but it also covers each neighbor of $v$ in $H$ a second time. Hence it is a perfect code precisely when $v$ is isolated in $H$, equivalently when $\deg_\Gamma(v)=1$. In that case the first assertion excludes cost zero, so the minimum is one.
\end{proof}

\begin{figure}[htbp]
\centering
\begin{tikzpicture}[x=8mm,y=8mm,font=\small,
  vertex/.style={circle,draw,minimum size=5mm,inner sep=0pt},
  code/.style={vertex,fill=black!75,text=white},
  failed/.style={densely dashed,thick}]
\node[code] (a) at (0,0) {$a$};
\node[vertex] (b) at (1,0) {$b$};
\node[vertex] (x) at (2,0) {$x$};
\node[code] (y) at (3,0) {$y$};
\draw (a)--(b) (x)--(y);
\draw[failed] (b)--(x);
\node at (1.5,-.7) {$r=0$: keep $C$};
\node[vertex] (w) at (5,0) {$w$};
\node[code] (u) at (6,0) {$u$};
\node[vertex] (v) at (7,0) {$v$};
\draw (w)--(u);
\draw[failed] (u)--(v);
\node at (6,-.7) {$r=1$: add $v$};
\end{tikzpicture}
\caption{Filled vertices form the initial code. Each dashed edge marks a single deletion.}
\label{fig:repair-cost-basics}
\end{figure}
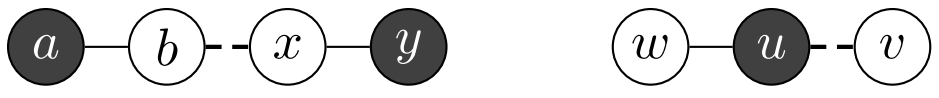

Figure~\ref{fig:repair-cost-basics} illustrates both cases of Lemma~\ref{l:elementary-repair}.

The graph-level extrema account for every edge deletion. On $P_4$, the unique perfect code is $\{a,y\}$. Deleting $ab$ gives the optimal repair $\{a,x\}$, and deleting $xy$ gives the symmetric repair $\{b,y\}$. The three edge costs in path order are therefore $2,0,2$, giving $R_{\min}(P_4)=R_{\max}(P_4)=2$. On $P_3$, either deletion isolates a leaf and costs one, so $R_{\min}(P_3)=R_{\max}(P_3)=1$.
In both examples, the minimum and maximum over perfect codes coincide because there is only one original perfect code. Graphs with several perfect codes may have different worst-case costs for different initial choices; Example~\ref{ex:depend} exhibits $R_{\min}(\Gamma)=2<R_{\max}(\Gamma)=4$.

\section{General Bounds, Parity, and Finiteness}
\begin{proposition}\label{p:gamma}
	Every perfect code $C$ of a graph $\Gamma$ satisfies $|C|=\gamma(\Gamma)$.
\end{proposition}
\begin{proof}
	The set $C$ dominates $\Gamma$, so $\gamma(\Gamma)\leq |C|$. Conversely, every dominating set $S$ must intersect $N_\Gamma[c]$ for each $c\in C$. These closed neighborhoods are pairwise disjoint because $C$ is a perfect code. Thus $|S|\geq |C|$, yielding equality.
\end{proof}
\begin{corollary}\label{c:equal-size}
	All perfect codes of a given graph have the same cardinality.
\end{corollary}
\begin{lemma}\label{l:gamma}
	For every finite simple graph $\Gamma$ and every edge $e\in E(\Gamma)$,
	\[
	\gamma(\Gamma)\leq\gamma(\Gamma-e)\leq\gamma(\Gamma)+1.
	\]
\end{lemma}
\begin{proof}
	Every dominating set of $\Gamma-e$ dominates $\Gamma$, giving the lower bound. For the upper bound, let $S$ be a minimum dominating set of $\Gamma$ and write $e=uv$. If $S$ fails to dominate $\Gamma-e$, only an endpoint of $e$ can have lost its dominator. Relabel the endpoints so that $v$ is undominated. Then $v\notin S$ and $u\in S$, so $S\cup\{v\}$ dominates $\Gamma-e$. In either case, $\gamma(\Gamma-e)\leq |S|+1$.
\end{proof}
\begin{lemma}\label{l:size}
	Let $e\in E(\Gamma)$, $C\in\PC(\Gamma)$, and $D\in\PC(\Gamma-e)$. Then $|D|\in\{|C|,|C|+1\}$.
\end{lemma}
\begin{proof}
	Apply Proposition~\ref{p:gamma} to both graphs and use Lemma~\ref{l:gamma}.
\end{proof}

\begin{proposition}\label{p:parity}
	If $r(\Gamma,C;e)<\infty$, then
	\[
	r(\Gamma,C;e)\equiv\gamma(\Gamma-e)-\gamma(\Gamma)\pmod 2.
	\]
\end{proposition}

%\begin{proof}
%	Choose an optimal repaired code $D$. By Proposition~\ref{p:gamma},
%	\[
%	r(\Gamma,C;e)=|C|+|D|-2|C\cap D|
%	=\gamma(\Gamma)+\gamma(\Gamma-e)-2|C\cap D|.
%	\]
%	Reducing modulo two proves the claim.
%\end{proof}

\begin{proof}
	Suppose that $r(\Gamma,C;e)<\infty$. Then there exists
	a perfect code $D$ of $\Gamma-e$ such that
	\[
	r(\Gamma,C;e)
	=|C\triangle D|
	=|C|+|D|-2|C\cap D|.
	\]
	By Proposition~\ref{p:gamma}, we have
	\[
	|C|=\gamma(\Gamma)
	\quad\text{and}\quad
	|D|=\gamma(\Gamma-e).
	\]
	Consequently,
	\[
	\begin{aligned}
		r(\Gamma,C;e)
		&=\gamma(\Gamma)+\gamma(\Gamma-e)-2|C\cap D|\\
		&\equiv \gamma(\Gamma)+\gamma(\Gamma-e)\pmod{2}\\
		&\equiv \gamma(\Gamma-e)-\gamma(\Gamma)\pmod{2},
	\end{aligned}
	\]
	as required.
\end{proof}

%\begin{corollary}
%	For any repairable edge $e$, the cost is even if and only if $\gamma(\Gamma-e)=\gamma(\Gamma)$, and odd if and only if $\gamma(\Gamma-e)=\gamma(\Gamma)+1$.
%\end{corollary}

Proposition~\ref{p:parity} establishes a connection between the parity
of a finite edge repair cost and the change in the domination
number under edge deletion.

\begin{corollary}\label{c:parity-cases}
	Let $C$ be a perfect code of a graph $\Gamma$ and let
	$e\in E(\Gamma)$. If $r(\Gamma,C;e)<\infty$, then
	\begin{enumerate}
		\item $r(\Gamma,C;e)$ is even if and only if
		$\gamma(\Gamma-e)=\gamma(\Gamma)$;
		
		\item $r(\Gamma,C;e)$ is odd if and only if
		$\gamma(\Gamma-e)=\gamma(\Gamma)+1$.
	\end{enumerate}
\end{corollary}

\begin{proof}
	By Lemma~\ref{l:gamma},
	\[
	0\leq \gamma(\Gamma-e)-\gamma(\Gamma)\leq 1.
	\]
	By Proposition~\ref{p:parity},
	\[
	r(\Gamma,C;e)
	\equiv \gamma(\Gamma-e)-\gamma(\Gamma)\pmod{2}.
	\]
	Since $\gamma(\Gamma-e)-\gamma(\Gamma)\in\{0,1\}$,
	it follows that $r(\Gamma,C;e)$ is even if and only if
	$\gamma(\Gamma-e)-\gamma(\Gamma)=0$, and odd if and only if
	$\gamma(\Gamma-e)-\gamma(\Gamma)=1$.
	This establishes both assertions.
\end{proof}

\begin{example}\label{ex:parity}
	In $K_3$ with code $\{a\}$, deleting $ab$ permits the repaired code $\{c\}$, giving cost two and leaving the domination number unchanged. In $P_3$ with its central vertex as code, deleting a pendant edge gives cost one and increases the domination number by one.
\end{example}
\begin{figure}[htbp]
	\centering
	\begin{tikzpicture}[scale=0.95]
		\node[code] (a) at (0,1.4) {$a$}; \node[vertex] (b) at (-.8,0) {$b$}; \node[vertex] (c) at (.8,0) {$c$};
		\draw[failed] (a)--(b); \draw[edge] (a)--(c)--(b);
		\draw[-{Stealth}] (1.3,.6)--(2.3,.6);
		\node[vertex] (a2) at (3.6,1.4) {$a$}; \node[vertex] (b2) at (2.8,0) {$b$}; \node[code] (c2) at (4.4,0) {$c$}; \draw[edge] (a2)--(c2)--(b2);
		\node at (1.8,-.65) {Even cost: $2$};
		\node[vertex] (x) at (0,-2) {$x$}; \node[code] (y) at (1,-2) {$y$}; \node[vertex] (z) at (2,-2) {$z$}; \draw[failed] (x)--(y); \draw[edge] (y)--(z);
		\draw[-{Stealth}] (2.5,-2)--(3.3,-2);
		\node[code] (x2) at (3.8,-2) {$x$}; \node[code] (y2) at (4.8,-2) {$y$}; \node[vertex] (z2) at (5.8,-2) {$z$}; \draw[edge] (y2)--(z2);
		\node at (2.9,-2.65) {Odd cost: $1$};
	\end{tikzpicture}
	\caption{The parity of a finite repair cost records whether the domination number increases. Filled vertices belong to the respective codes; dashed edges are deleted.}
\end{figure}
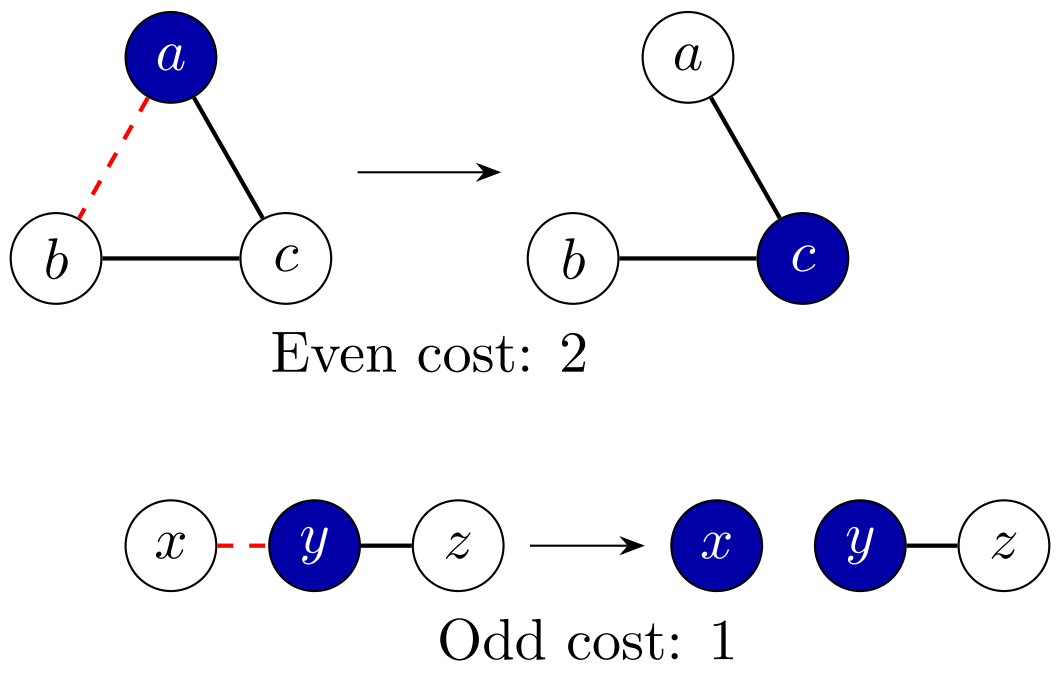

\begin{proposition}\label{p:bounds}
	If $R(\Gamma,C)<\infty$, then
	\[
	1\leq R(\Gamma,C)\leq\min\{n-1,2\gamma(\Gamma)+1\}.
	\]
\end{proposition}
\begin{proof}
	There is at least one edge incident with $C$, so $R(\Gamma,C)\geq1$. For any edge $e$ and optimal repaired code $D$, Lemma~\ref{l:size} yields
	\[
	|C\dist D|\leq |C|+|D|\leq2\gamma(\Gamma)+1.
	\]
	Suppose instead that $r(\Gamma,C;e)=n$. Then $D=V(\Gamma)\setminus C$. Since the cost is positive, write $e=uv$ with $u\in C$. In $\Gamma-e$, the vertex $v$ has no code vertex of $C$ in its closed neighborhood. Thus $N_{\Gamma-e}[v]\subseteq D$. As $D$ is a perfect code, $N_{\Gamma-e}[v]=\{v\}$, so $v$ is a leaf in $\Gamma$. Lemma~\ref{l:elementary-repair} implies $r(\Gamma,C;e)=1$, contradicting $n\geq2$. Therefore $r(\Gamma,C;e)\leq n-1$ for all $e$.
\end{proof}
\begin{theorem}\label{t:pendant-cycle}
For every $q\ge1$, let $\Gamma_q$ be the cycle $v_1v_2\cdots v_{3q}v_1$ with a leaf $w$ adjacent to $v_1$. Then $C=\{v_1,v_4,\ldots,v_{3q-2}\}$ is a perfect code, $\gamma(\Gamma_q)=q$, and $R(\Gamma_q,C)=2q+1$.
\end{theorem}
\begin{proof}
Each cycle vertex is covered exactly once by $C$, and $v_1$ covers $w$. Thus $C$ is a perfect code and Proposition~\ref{p:gamma} gives $\gamma(\Gamma_q)=q$.

First we verify finiteness. Deleting $wv_1$ permits the code $C\cup\{w\}$, and deleting a cycle edge disjoint from $C$ preserves $C$. If a cycle edge incident with $C$ is deleted, the remaining cycle edges form a path. Its two arms measured from $v_1$ have lengths congruent to $0$ and $2$ modulo three, in some order (an arm may be empty). Indeed, for the edge $v_iv_{i+1}$ with cyclic indices the lengths are $i-1$ and $3q-i$, and incidence with $C$ means $i\equiv0$ or $1\pmod3$. Select $w$ and, on each arm, the vertices at distances $2,5,8,\ldots$ from $v_1$. These centers cover consecutive triples along the arm, with a final pair when its length is $2$ modulo three. The vertex $w$ covers itself and $v_1$, so this is a perfect code after deletion.

Finally, deleting $v_{3q}v_1$ leaves the path in the order $w,v_1,\ldots,v_{3q}$. By Lemma~\ref{l:pathcodes}, its unique perfect code is $\{w,v_3,v_6,\ldots,v_{3q}\}$. This code is disjoint from $C$ and has size $q+1$, giving cost $2q+1$. Proposition~\ref{p:bounds} gives the matching upper bound.
\end{proof}

\begin{remark}\label{rem:bounds}
Proposition~\ref{p:bounds} is a general upper estimate; it need not be attained for a given graph. For $P_4$, the unique perfect code has maximum repair cost $2$, whereas
$\min\{n-1,2\gamma(P_4)+1\}=\min\{3,5\}=3$.
The bound is attained for some graphs, including $K_2$ and $K_3$.
Theorem~\ref{t:pendant-cycle} attains the term $2\gamma(\Gamma)+1$ for every positive domination number. Characterizing all equality cases remains open.
Likewise, the spread bounds below are upper estimates; no general equality characterization is asserted.
\end{remark}

\begin{corollary}\label{c:even-maximum}
	Let $C$ be a perfect code of a graph $\Gamma$ with
	$E(\Gamma)\neq\emptyset$. If every edge is repairable and
	no edge deletion increases $\gamma(\Gamma)$, then
	$R(\Gamma,C)$ is even and
	\[
	R(\Gamma,C)\leq 2\gamma(\Gamma).
	\]
\end{corollary}

\begin{proof}
	Let $e\in E(\Gamma)$ be arbitrary. Since every edge is
	repairable, there exists a perfect code
	$D\in\PC(\Gamma-e)$ such that
	\[
	r(\Gamma,C;e)=|C\triangle D|.
	\]
	By Lemma~\ref{l:size}, we have
	\[
	|D|\in\{|C|,|C|+1\}.
	\]
	Moreover, since no edge deletion increases the domination
	number, we have
	\[
	\gamma(\Gamma-e)=\gamma(\Gamma).
	\]
	It follows from Proposition~\ref{p:gamma} that
	\[
	|D|=\gamma(\Gamma-e)=\gamma(\Gamma)=|C|.
	\]
	Hence,
	\[
	\begin{aligned}
		r(\Gamma,C;e)
		&=|C|+|D|-2|C\cap D|\\
		&=2|C|-2|C\cap D|\\
		&\leq 2|C|\\
		&=2\gamma(\Gamma).
	\end{aligned}
	\]
	Furthermore, by Proposition~\ref{p:parity},
	\[
	r(\Gamma,C;e)
	\equiv \gamma(\Gamma-e)-\gamma(\Gamma)
	\equiv 0\pmod{2}.
	\]
	Thus, $r(\Gamma,C;e)$ is an even nonnegative integer
	not exceeding $2\gamma(\Gamma)$.
	
	Since $e$ was chosen arbitrarily, the same conclusion
	holds for every edge of $\Gamma$. Therefore,
	\[
	R(\Gamma,C)
	=\max_{e\in E(\Gamma)}r(\Gamma,C;e)
	\leq 2\gamma(\Gamma).
	\]
	Moreover, being the maximum of a nonempty finite set of
	even integers, $R(\Gamma,C)$ is even.
\end{proof}

\section[Structural Characterizations and Dependence on the Initial Code]{Structural Characterizations and\\Dependence on the Initial Code}
The following result shows that the finiteness of the maximum
edge repair cost is independent of the choice of a perfect code.

\begin{proposition}\label{p:finite}
	Let $\Gamma$ be a graph with $\PC(\Gamma)\neq\varnothing$
	and $E(\Gamma)\neq\varnothing$. For every $C\in\PC(\Gamma)$,
	\[
	R(\Gamma,C)=\infty
	\quad\Longleftrightarrow\quad
	\exists e\in E(\Gamma)\text{ such that }
	\PC(\Gamma-e)=\varnothing.
	\]
	Consequently, the following statements are equivalent:
	\begin{enumerate}
		\item $R_{\min}(\Gamma)<\infty$;
		\item $R_{\max}(\Gamma)<\infty$;
		\item $\PC(\Gamma-e)\neq\varnothing$ for every
		$e\in E(\Gamma)$.
	\end{enumerate}
\end{proposition}

\begin{proof}
	Let $C\in\PC(\Gamma)$. By definition,
	\[
	R(\Gamma,C)
	=\max_{e\in E(\Gamma)}r(\Gamma,C;e).
	\]
	Since $E(\Gamma)$ is finite and nonempty, we have
	$R(\Gamma,C)=\infty$ if and only if
	$r(\Gamma,C;e)=\infty$ for some $e\in E(\Gamma)$.
	
	For any edge $e\in E(\Gamma)$, if
	$\PC(\Gamma-e)\neq\varnothing$, then there exists a
	perfect code $D\in\PC(\Gamma-e)$ such that
	\[
	r(\Gamma,C;e)
	\leq |C\triangle D|
	\leq |V(\Gamma)|<\infty.
	\]
	On the other hand, if $\PC(\Gamma-e)=\varnothing$,
	then $r(\Gamma,C;e)=\infty$ by definition. Hence,
	\[
	R(\Gamma,C)=\infty
	\quad\Longleftrightarrow\quad
	\exists e\in E(\Gamma)\text{ such that }
	\PC(\Gamma-e)=\varnothing.
	\]
	
	We now establish the equivalence of the three statements.
	
	Suppose that {\rm (1)} holds. Then there exists a perfect
	code $C\in\PC(\Gamma)$ such that $R(\Gamma,C)<\infty$.
	By the criterion for $R(\Gamma,C)=\infty$ just proved, $\PC(\Gamma-e)\neq\varnothing$
	for every $e\in E(\Gamma)$. Thus, {\rm (3)} holds.
	
	Suppose that {\rm (3)} holds. By the same criterion, $R(\Gamma,C)<\infty$ for every
	$C\in\PC(\Gamma)$. Since $\PC(\Gamma)$ is finite and nonempty,
	\[
	R_{\max}(\Gamma)
	=\max_{C\in\PC(\Gamma)}R(\Gamma,C)<\infty.
	\]
	Hence, {\rm (2)} holds.
	
	Finally, {\rm (2)} implies {\rm (1)}, since
	\[
	R_{\min}(\Gamma)\leq R_{\max}(\Gamma).
	\]
	This completes the proof.
\end{proof}

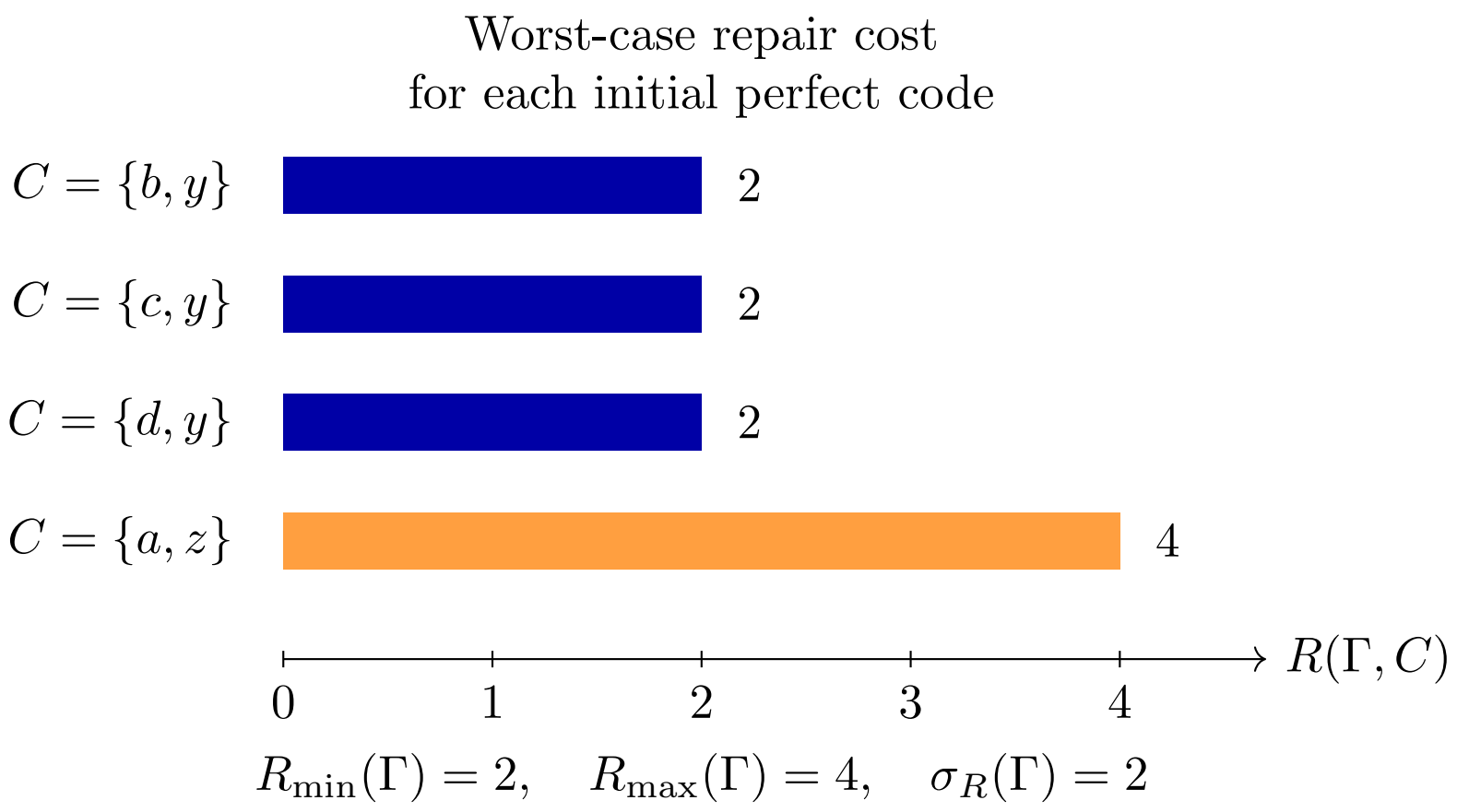
\begin{figure}[htbp]
\centering
\begin{tikzpicture}[x=1.5cm,y=.85cm,font=\small]
\draw[->] (0,0)--(4.7,0) node[right] {$R(\Gamma,C)$};
\foreach \k in {0,1,2,3,4}{
  \draw (\k,.07)--(\k,-.07) node[below] {$\k$};
}
\foreach \row/\lab/\cost in {4/{\{b,y\}}/2,3/{\{c,y\}}/2,2/{\{d,y\}}/2,1/{\{a,z\}}/4}{
  \node[anchor=east] at (-.15,\row) {$C=\lab$};
  \ifnum\cost=4
    \fill[orange!75] (0,\row-.24) rectangle (\cost,\row+.24);
  \else
    \fill[blue!65!black] (0,\row-.24) rectangle (\cost,\row+.24);
  \fi
  \node[anchor=west] at (\cost+.08,\row) {$\cost$};
}
\node[align=center] at (2,5) {Worst-case repair cost\\for each initial perfect code};
\node[align=center] at (2,-1.0) {$R_{\min}(\Gamma)=2,\quad R_{\max}(\Gamma)=4,\quad \sigma_R(\Gamma)=2$};
\end{tikzpicture}
\caption{Dependence of the repair budget on the initial perfect code. Blue bars mark the three minimizing codes; the orange bar marks the maximizing code. Each code has cardinality two, but the maximum repair costs differ.}
\label{fig:code-choice-costs}
\end{figure}

\begin{corollary}\label{c:extrema-bounds}
	If the parameters are finite, Proposition~\ref{p:bounds} gives
	\[
	1\leq R_{\min}(\Gamma)\leq R_{\max}(\Gamma)\leq\min\{n-1,2\gamma(\Gamma)+1\}.
	\]
\end{corollary}
\begin{proposition}\label{p:star-forest}
	The following are equivalent: $R_{\min}(\Gamma)=1$; $R_{\max}(\Gamma)=1$; and $\Gamma$ is a star forest.
\end{proposition}
\begin{proof}
	If $R(\Gamma,C)=1$, every vertex outside $C$ must be a leaf by Lemma~\ref{l:elementary-repair}. Since $C$ is independent and every vertex outside $C$ has exactly one neighbor in $C$, each component is a star centered at a vertex of $C$ (including isolated vertices). Conversely, a star with at least three vertices has its center as its unique perfect code, an isolated vertex must belong to every code, and either endpoint of $K_2$ is a center. Thus every edge of a star forest joins a code vertex to a leaf outside the code, and the unit-cost criterion gives maximum repair cost one.
\end{proof}
\begin{example}[Dependence on the original code]\label{ex:depend}
	Let $\Gamma$ consist of a complete graph on $\{a,b,c,d\}$ with the path $a x y z$ attached at $a$. Its perfect codes are exactly
	\[
	\{b,y\},\quad\{c,y\},\quad\{d,y\},\quad\{a,z\}.
	\]
	Indeed, exactly one clique vertex must be selected to dominate $b,c,d$. If $a$ is selected, $x$ is covered and $z$ must be selected. Otherwise $x$ must be covered by $y$, since selecting $x$ would dominate $a$ twice. These choices give precisely the four displayed codes.
	
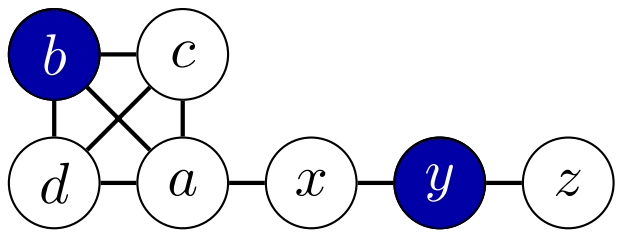
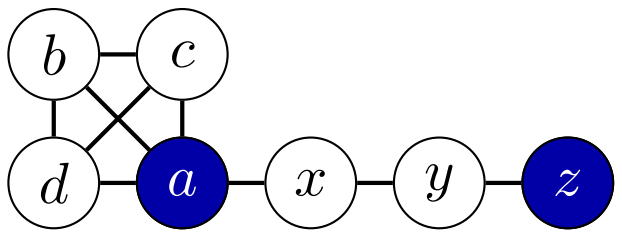
\begin{center}
		\begin{tikzpicture}[scale=.85]
			\foreach \s/\lab in {0/1,7/2}{
				\begin{scope}[xshift=\s cm]
					\node[vertex] (a\lab) at (1,0) {$a$}; \node[vertex] (b\lab) at (0,1) {$b$}; \node[vertex] (c\lab) at (1,1) {$c$}; \node[vertex] (d\lab) at (0,0) {$d$};
					\node[vertex] (x\lab) at (2,0) {$x$}; \node[vertex] (y\lab) at (3,0) {$y$}; \node[vertex] (z\lab) at (4,0) {$z$};
					\draw[edge] (a\lab)--(b\lab)--(c\lab)--(d\lab)--(a\lab) (a\lab)--(c\lab) (b\lab)--(d\lab) (a\lab)--(x\lab)--(y\lab)--(z\lab);
					\ifnum\lab=1 \node[code] at (b\lab) {$b$}; \node[code] at (y\lab) {$y$}; \node at (2,-.7) {$C_1=\{b,y\}$};
					\else \node[code] at (a\lab) {$a$}; \node[code] at (z\lab) {$z$}; \node at (2,-.7) {$C_2=\{a,z\}$}; \fi
				\end{scope}
			}
		\end{tikzpicture}
	\end{center}
	The following table gives the exact costs for the two representative codes. The other two codes are equivalent to $C_1$ under a permutation of $b,c,d$.
	\begin{center}
		\begin{tabular}{ccccc}
			\toprule
			Deleted edge & Optimal code for $C_1$ & Cost & Optimal code for $C_2$ & Cost\\
			\midrule
			$ab$ & $\{c,y\}$ & 2 & $\{c,y\}$ & 4\\
			$ac$ & $\{b,y\}$ & 0 & $\{b,y\}$ & 4\\
			$ad$ & $\{b,y\}$ & 0 & $\{b,y\}$ & 4\\
			$bc$ & $\{d,y\}$ & 2 & $\{a,z\}$ & 0\\
			$bd$ & $\{c,y\}$ & 2 & $\{a,z\}$ & 0\\
			$cd$ & $\{b,y\}$ & 0 & $\{a,z\}$ & 0\\
			$ax$ & $\{b,y\}$ & 0 & $\{a,y\}$ & 2\\
			$xy$ & $\{a,y\}$ & 2 & $\{a,z\}$ & 0\\
			$yz$ & $\{b,y,z\}$ & 1 & $\{b,y,z\}$ & 3\\
			\bottomrule
		\end{tabular}
	\end{center}
	These minima can be checked by enumerating subsets of the seven vertices and testing the closed-neighborhood condition after each deletion. Hence $R(\Gamma,C_1)=2$ and $R(\Gamma,C_2)=4$, so
	\[
	R_{\min}(\Gamma)=2,\qquad R_{\max}(\Gamma)=4.
	\]
\end{example}
The four initial choices are compared in Figure~\ref{fig:code-choice-costs}. The labels agree with Example~\ref{ex:depend}; they correspond to $(a,b,c,d,x,y,z)=(v_1,v_2,v_3,v_4,v_6,v_5,v_7)$ in the numerical labeling. The graph drawings in that example show a best and a worst initial choice, while the bars below compare all four choices.

The example establishes that selecting a minimum-cardinality dominating set which is also a perfect code need not minimize the subsequent repair budget.

\begin{proposition}\label{p:cardinality-versus-repair}
There exists a connected, fully edge-repairable graph $\Gamma$ whose perfect codes all have cardinality $\gamma(\Gamma)$ but for which
\[
R_{\min}(\Gamma)<R_{\max}(\Gamma).
\]
In particular, minimum cardinality does not imply minimum worst-case repair cost.
\end{proposition}
\begin{proof}
Take the graph in Example~\ref{ex:depend}. Its four perfect codes each have cardinality two, so Proposition~\ref{p:gamma} gives $\gamma(\Gamma)=2$. The edge-cost table and the symmetry among $b,c,d$ show that all four codes have finite maximum repair costs, with values $2,2,2,4$. Hence $\Gamma$ is fully edge-repairable and
\[
R_{\min}(\Gamma)=2<4=R_{\max}(\Gamma).
\]
The code $\{a,z\}$ therefore has minimum cardinality but does not minimize the worst-case repair cost.
\end{proof}

For this graph, selecting any of $\{b,y\}$, $\{c,y\}$, or $\{d,y\}$ instead of $\{a,z\}$ reduces the worst-case budget by $(4-2)/4=50\%$. This comparison concerns the maximum cost over a single deleted edge, rather than an average over failures.

\subsection{Smallest order with unequal repair costs}

\begin{theorem}\label{t:small-order}
Among finite simple graphs with at least one edge and a perfect code, the least order for which $R_{\min}(\Gamma)<R_{\max}(\Gamma)$ is seven. Up to isomorphism, exactly two graphs of order seven have this property:
\begin{itemize}
\item $K_4$ with a path of length three attached at one endpoint, with extrema $(2,4)$;
\item $K_5-ab$ with a path of length two attached at $b$, with extrema $(3,4)$.
\end{itemize}
\end{theorem}
\begin{proof}[Proof by exhaustive enumeration and orbit counting]
Appendix~\ref{app:census} gives the complete program. On the fixed vertex set $\{0,\ldots,n-1\}$ it tests all $2^{\binom n2}-1$ nonempty edge sets, including disconnected graphs. For each graph with at least two perfect codes it computes the two extrema directly. A deletion admitting no perfect code gives every initial code infinite worst-case cost. The numbers of labelled graphs with unequal extrema are
\[
\begin{array}{c|rrrrrr}
n&2&3&4&5&6&7\\\hline
\text{Graphs tested}&1&7&63&1023&32767&2097151\\
(m,R_{\min},R_{\max})=(9,2,4)&0&0&0&0&0&840\\
(m,R_{\min},R_{\max})=(11,3,4)&0&0&0&0&0&840\\
\text{Other triples with unequal extrema}&0&0&0&0&0&0
\end{array}
\]
The two examples in this section realize the indicated triples. In the first graph every automorphism fixes the unique degree-four vertex and the attached path pointwise, and freely permutes the remaining three clique vertices. In the second, the leaf $y$, its neighbor $x$, and the other neighbor $b$ of $x$ are fixed. The vertex $a$ is the only degree-three vertex, and the remaining three vertices are freely permutable. Each graph therefore has $3!$ automorphisms and accounts for $7!/3!=840$ distinct labellings. Their edge counts distinguish the two isomorphism classes, and their labellings exhaust the census.
\end{proof}

A second graph of minimum order realizes spread one. On $\{a,b,c,d,e\}$ take $K_5-ab$ and attach the path $bxy$. Its four perfect codes, in the column order used below, are $C_1=\{c,y\}$, $C_2=\{d,y\}$, $C_3=\{e,y\}$, and $C_4=\{a,x\}$. Applying the same closed-neighborhood test to all $128$ vertex subsets before and after each deletion gives the following costs.
\[
\begin{array}{c|rrrr}
f&C_1&C_2&C_3&C_4\\\hline
ac & 2 & 0 & 0 & 4 \\
ad & 0 & 2 & 0 & 4 \\
ae & 0 & 0 & 2 & 4 \\
bc & 2 & 0 & 0 & 0 \\
bd & 0 & 2 & 0 & 0 \\
be & 0 & 0 & 2 & 0 \\
cd & 2 & 2 & 0 & 0 \\
ce & 2 & 0 & 2 & 0 \\
de & 0 & 2 & 2 & 0 \\
bx & 0 & 0 & 0 & 2 \\
xy & 3 & 3 & 3 & 1 \\
\hline
\text{Maximum}&3&3&3&4
\end{array}
\]
\begin{figure}[htbp]
\centering
\begin{tikzpicture}
\foreach \shift/\stage in {0/1,5/2}{
\begin{scope}[xshift=\shift cm]
\foreach \v/\i/\j in {a/0/1,b/2/1,c/1/2,d/0/0,e/2/0,x/3/1,y/4/1}{
\node[vertex] (\v\stage) at (\i,\j) {$\v$};}
\foreach \u/\v in {a/d,a/e,b/c,b/d,b/e,c/d,c/e,d/e,b/x,x/y}{\draw (\u\stage)--(\v\stage);}
\ifnum\stage=1
\draw[failed] (a1)--(c1);
\node[code] at (a1) {$a$}; \node[code] at (x1) {$x$};
\node at (2,-.7) {$\Gamma,\quad C_4=\{a,x\}$};
\else
\node[code] at (d2) {$d$}; \node[code] at (y2) {$y$};
\node at (2,-.7) {$\Gamma-ac,\quad D=\{d,y\}$};
\fi
\end{scope}}
\end{tikzpicture}
\caption{Deleting the dashed edge $ac$ forces both vertices of $C_4$ to be replaced. Filled vertices indicate the selected codes.}
\label{fig:seven-one}
\end{figure}
Indeed, the only perfect codes after deleting $ac$ are $\{d,y\}$ and $\{e,y\}$. Either repair retains $y$ from $C_1$ and retains neither vertex of $C_4$, giving respective costs two and four. Thus this eleven-edge graph has $(R_{\min},R_{\max},\sigma_R)=(3,4,1)$, whereas the nine-edge graph of Example~\ref{ex:depend} has $(2,4,2)$.

\subsection{Stability and Repair-Cost Spread}
\begin{proposition}\label{p:stability}
	If $C_1,C_2\in\PC(\Gamma)$ have finite maximum repair costs, then
	\[
	|R(\Gamma,C_1)-R(\Gamma,C_2)|\leq |C_1\dist C_2|.
	\]
\end{proposition}
\begin{proof}
	For every edge $e$, let $D$ be an optimal repair of $C_2$. The triangle inequality yields
	\[
	r(\Gamma,C_1;e)\leq |C_1\dist D|\leq |C_1\dist C_2|+r(\Gamma,C_2;e).
	\]
	Take the maximum over $e$ and then interchange $C_1$ and $C_2$.
\end{proof}
\begin{definition}\label{def:spread}
	For a graph with finite maximum repair costs, define the \emph{repair-cost spread} by
	\[
	\sigma_R(\Gamma)=R_{\max}(\Gamma)-R_{\min}(\Gamma).
	\]
\end{definition}
For the nonempty family of perfect codes, define its symmetric-difference diameter by
\[
\operatorname{diam}_{\triangle}\PC(\Gamma)
=\max_{C_1,C_2\in\PC(\Gamma)}|C_1\triangle C_2|.
\]
This diameter is zero when the perfect code is unique. The stability inequality yields a bound that also accounts for the distances between the available codes.

\begin{corollary}\label{c:spread-diameter}
If $\Gamma$ has finite maximum repair costs, then
\[
0\leq\sigma_R(\Gamma)\leq
\operatorname{diam}_{\triangle}\PC(\Gamma)\leq2\gamma(\Gamma).
\]
Moreover, $\sigma_R(\Gamma)=0$ precisely when the maximum repair cost is independent of the initial perfect code.
\end{corollary}
\begin{proof}
Choose $C_+$ and $C_-$ attaining $R_{\max}(\Gamma)$ and $R_{\min}(\Gamma)$, respectively. Proposition~\ref{p:stability} gives
\[
\sigma_R(\Gamma)\leq |C_+\triangle C_-|
\leq\operatorname{diam}_{\triangle}\PC(\Gamma).
\]
Every perfect code has cardinality $\gamma(\Gamma)$, so
$|C_1\triangle C_2|\leq |C_1|+|C_2|=2\gamma(\Gamma)$.
The last assertion follows from the definition of the spread.
\end{proof}

The order bound improves when the graph is not a star forest, since no initial code can then have maximum repair cost one.

\begin{corollary}\label{c:spread-order}
Suppose that $\Gamma$ has order $n$, at least one edge, and finite maximum repair costs. If $\Gamma$ is a star forest, then $\sigma_R(\Gamma)=0$. Otherwise,
\[
0\leq\sigma_R(\Gamma)\leq
\min\{n-3,\operatorname{diam}_{\triangle}\PC(\Gamma)\}
\leq\min\{n-3,2\gamma(\Gamma)\}.
\]
\end{corollary}
\begin{proof}
Every perfect code of a star forest has maximum repair cost one, giving spread zero.
If $\Gamma$ is not a star forest, the unit-cost characterization and integrality imply $R_{\min}(\Gamma)\geq2$.
Proposition~\ref{p:bounds} gives $R_{\max}(\Gamma)\leq n-1$.
Consequently,
$\sigma_R(\Gamma)=R_{\max}(\Gamma)-R_{\min}(\Gamma)\leq n-3$.
Combine this with Corollary~\ref{c:spread-diameter}.
\end{proof}
Example~\ref{ex:depend} has $\sigma_R(\Gamma)=2$.

\subsection{Repair by partitions of closed neighborhoods}
Fix $C\in\PC(\Gamma)$ and an edge $e=uv$ with $u\in C$, and write $H=\Gamma-e$. For $W\subseteq V(H)$, let $\nu_H(W)$ be the minimum size of a set $B\subseteq W$ such that the closed neighborhoods $N_H[b]$, indexed by $b\in B$, are pairwise disjoint and have union $W$. Set $\nu_H(\varnothing)=0$ and take the minimum to be $\infty$ when no such partition exists. For $A\subseteq C$, put
\[
W_A=\{v\}\cup\bigcup_{a\in A}N_H[a].
\]

\begin{theorem}\label{t:repair-partition}
Under these hypotheses,
\[
r(\Gamma,C;e)=\min_{A\subseteq C}\bigl(|A|+\nu_H(W_A)\bigr).
\]
\end{theorem}
\begin{proof}
The neighborhoods centered at $C$ partition $V(H)\setminus\{v\}$. Hence those centered at $C\setminus A$ partition $V(H)\setminus W_A$. For any perfect code $D$ of $H$, take $A=C\setminus D$. The centers in $D\setminus C$ must then partition $W_A$ by their closed neighborhoods. Consequently $|A|+\nu_H(W_A)\le |C\dist D|$.

Conversely, suppose $B$ gives a partition of $W_A$. The two complementary partitions show that $D=(C\setminus A)\cup B$ is a perfect code of $H$. Since $B\subseteq W_A$ is disjoint from $C\setminus A$, its actual repair cost is
\[
|C\dist D|=|A|+|B|-2|A\cap B|\le |A|+|B|.
\]
Minimizing gives the reverse inequality. The constructions also show that a perfect code of $H$ exists exactly when at least one term is finite.
\end{proof}

The formula gives an exact algorithm. For each $A$, compute $\nu_H(W_A)$ recursively. On a nonempty remaining set $U$, choose any $x\in U$ and branch over centers $b\in N_H[x]$ with $N_H[b]\subseteq U$. Each branch contributes $1+\nu_H(U\setminus N_H[b])$. Take the minimum, using $\infty$ if there is no branch. Every partition covers $x$ in exactly one block, and every branch removes at least one vertex. Values for repeated sets $U$ can be stored and reused across choices of $A$.

For $A=\varnothing$, the uncovered set is $\{v\}$, which admits a partition precisely when $v$ is isolated in $H$. This recovers the unit-cost criterion in Lemma~\ref{l:elementary-repair}. More generally, if a repair exists, put $\delta=\gamma(H)-\gamma(\Gamma)\in\{0,1\}$ by Lemma~\ref{l:gamma}. For a repaired code $D$, take its actual removals $A=C\setminus D$ and additions $B=D\setminus C$. Proposition~\ref{p:gamma} gives
\[
|B|-|A|=\delta,\qquad |C\dist D|=2|A|+\delta,
\]
recovering Proposition~\ref{p:parity}.

In particular, cost three occurs exactly when $v$ is not a leaf in $\Gamma$ and some singleton $A\subseteq C$ admits a partition of $W_A$ by two centers outside $C$. Such a partition forces $\delta=1$, and the leaf criterion excludes the only smaller odd cost. Cost four occurs exactly when some two-element $A\subseteq C$ admits a partition by two centers outside $C$ and no singleton $A$ admits a partition by one center outside $C$. Here $\delta=0$, incidence with $C$ excludes zero, and the singleton test excludes cost two. In both cases necessity follows by taking the actual removals and additions of an optimal repair.

\subsection{Structural Characterization of Repair Cost Two}
\begin{theorem}\label{t:two}
	Let $e=uv$ with $u\in C$ and $\deg_\Gamma(v)\geq2$, and put $H=\Gamma-e$. Then $r(\Gamma,C;e)=2$ if and only if there exist $a\in C$ and $b\notin C$ such that
	\[
	N_H[b]=N_H[a]\cup\{v\}.
	\]
\end{theorem}
\begin{proof}
	A repair of cost two has one removal and one addition by the identity $|C\dist D|=2|C\setminus D|+\delta$ above. With $A=\{a\}$, Theorem~\ref{t:repair-partition} and its proof give the required single block $N_H[b]=W_A$ with $b\notin C$.
Conversely, if the condition holds, the replacement produces a perfect code of $H$. Its cost is two, and by Lemma~\ref{l:elementary-repair} the assumption that $v$ is not a leaf excludes cost one. Cost zero is impossible because $u\in C$.
\end{proof}
\begin{corollary}\label{c:Rtwo}
	The equality $R(\Gamma,C)=2$ holds if and only if both of the following conditions hold:
	\begin{enumerate}
		\item At least one vertex outside $C$ has degree at least two.
		\item For every $e=uv$ with $u\in C$ and $\deg_\Gamma(v)\geq2$, there exist $a\in C$ and $b\notin C$ satisfying the condition in Theorem~\ref{t:two}.
	\end{enumerate}
\end{corollary}
\begin{proof}
	By Lemma~\ref{l:elementary-repair}, edges disjoint from $C$ have cost zero, and edges joining code vertices to leaves outside $C$ have cost one. Theorem~\ref{t:two} characterizes the remaining edges of cost two. Condition 1 ensures that at least one such edge exists.
\end{proof}

\subsection{Irreparable Edges}
Define $B(\Gamma)=\{e\in E(\Gamma):\PC(\Gamma-e)=\varnothing\}$ and $E_C(\Gamma)=\{uv\in E(\Gamma):\{u,v\}\cap C\neq\varnothing\}$.
\begin{proposition}\label{p:bad-edges}
	\[
	B(\Gamma)\subseteq \bigcap_{C\in\PC(\Gamma)} E_C(\Gamma).
	\]
\end{proposition}
\begin{proof}
	If an edge has neither endpoint in some perfect code, that code survives its deletion. Such an edge cannot belong to $B(\Gamma)$.
\end{proof}
\begin{corollary}\label{c:edge-avoidance}
	If $\bigcap_{C\in\PC(\Gamma)}E_C(\Gamma)=\varnothing$, then $R_{\max}(\Gamma)<\infty$.
\end{corollary}
\begin{remark}\label{rem:edge-avoidance}
	The converse need not hold: $K_2$ has a single edge incident with every perfect code, but deletion of that edge is repairable.
\end{remark}

\section{Exact Results for Graph Families}
\subsection{Complete graphs}
\begin{proposition}\label{p:complete}
	For every perfect code $C$ of $K_n$, $n\geq2$,
	\[
	R(K_n,C)=\begin{cases}1,&n=2,\\2,&n\geq3.\end{cases}
	\]
\end{proposition}
\begin{proof}
	Every perfect code is a singleton, say $C=\{u\}$. For $n=2$, deleting the only edge isolates both vertices, so the repair cost is one. For $n\geq3$, deleting an edge not incident with $u$ preserves $C$. If $e=uv$, any vertex $w$ distinct from $u,v$ remains universal in $K_n-e$, so $\{w\}$ is a repaired code of cost two. Lemma~\ref{l:elementary-repair} excludes cost one because $v$ is not a leaf. The maximum is therefore two.
\end{proof}
\begin{corollary}\label{c:complete}
	For $n\geq3$, $R_{\min}(K_n)=R_{\max}(K_n)=2$.
\end{corollary}

\subsection{Paths}
\begin{lemma}\label{l:pathcodes}
	Let $P_n=v_1\cdots v_n$. Every perfect code consists of all vertices whose indices are congruent to one fixed residue $s\in\{1,2\}$ modulo three, with the following allowable residues:
	\[
	\begin{array}{c|c}
		n\pmod3 & s\\\hline
		0 & 2\\
		1 & 1\\
		2 & 1,2
	\end{array}
	\]
	This holds for every $n\geq1$.
\end{lemma}
\begin{proof}
	The first code vertex is $v_1$ or $v_2$ (when the latter exists), since $v_1$ must be dominated. Successive code vertices must be at distance exactly three: a smaller distance gives overlapping closed neighborhoods, and a larger distance leaves a vertex undominated. The last code vertex must be $v_{n-1}$ or $v_n$. These endpoint conditions give the displayed residues, and each allowed choice produces a perfect code. For $n=1$ the only code is $\{v_1\}$.
\end{proof}
\begin{theorem}\label{t:path}
	For $n\geq2$ and every perfect code $C$ of $P_n$,
	\[
	R(P_n,C)=\left\lfloor\frac{2n-1}{3}\right\rfloor
		=\begin{cases}2q-1,&n=3q,\\2q,&n=3q+1,\\2q+1,&n=3q+2.\end{cases}
	\]
\end{theorem}
\begin{proof}
	Deleting $e_i=v_iv_{i+1}$ produces paths on $i$ and $n-i$ vertices. A minimum repair is obtained by minimizing the symmetric difference independently on each component, because the components have disjoint vertex sets. Both admit perfect codes by Lemma~\ref{l:pathcodes}.
	
	Write $i=3t+j$ with $j\in\{0,1,2\}$. The costs are given in Table~\ref{tab:pathcost}; only indices $1\leq i<n$ are used. For $n=3q+2$, take the original code with residue one; reflection treats the other code.
	\begin{table}[htbp]
		\centering
		\caption{Exact single-edge repair costs of paths.}
		\label{tab:pathcost}
		\begin{tabular}{cccc}
			\toprule
			Order & $i=3t$ & $i=3t+1$ & $i=3t+2$\\
			\midrule
			$3q$ & $0$ & $2t+1$ & $2(q-t)-1$\\
			$3q+1$ & $2t$ & $2(q-t)$ & $0$\\
			$3q+2$ & $2t$ & $2(q-t)+1$ & $0$\\
			\bottomrule
		\end{tabular}
	\end{table}
	To justify the table, first suppose $n=3q$. The original code has residue two. For $i=3t$, neither endpoint is a code vertex, so the cost is zero. For $i=3t+1$, the left component has a unique perfect code, of size $t+1$ and disjoint from its original $t$ code vertices. The right component retains its original code. Thus the cost is $2t+1$. For $i=3t+2$, reflection gives the value $2(q-t)-1$.
	
	Now suppose $n=3q+1$ or $n=3q+2$, with original residue one. When $i=3t+2$ the edge avoids the code. When $i=3t$, the left component has a unique perfect code of size $t$, disjoint from the $t$ original code vertices; the right component retains its original code. The cost is $2t$. When $i=3t+1$, the left component retains its original code. The right component has a unique perfect code disjoint from the original code there. For $n=3q+1$ both right-component codes have size $q-t$; for $n=3q+2$ their sizes are $q-t$ and $q-t+1$. These give the two remaining entries.
	
	Taking the maximum over valid indices yields $2q-1$, $2q$, and $2q+1$, respectively. For $n=3q$, the maximum occurs at $i=2$ when $q\geq1$ (including $n=3$); for $n=3q+1$ or $3q+2$ it occurs at $i=1$.
\end{proof}
\begin{corollary}\label{c:path}
	For $n\geq2$, $R_{\min}(P_n)=R_{\max}(P_n)=\lfloor(2n-1)/3\rfloor$.
\end{corollary}
\begin{example}\label{ex:path-six}
	The unique perfect code of $P_6$ is $\{v_2,v_5\}$. The costs in path order are $1,3,0,3,1$, so $R(P_6,C)=3$. Deleting $v_2v_3$ gives the optimal repair $\{v_2,v_3,v_6\}$.
\end{example}

\subsection{Cycles}
\begin{theorem}\label{t:cycle}
	A cycle $C_n$ admits a perfect code if and only if $n=3q$ for an integer $q\geq1$. For every such code $C$,
	\[
	R(C_{3q},C)=2q.
	\]
\end{theorem}
\begin{proof}
	Every closed neighborhood of a cycle has size three, so $3|C|=n$. Conversely, choosing every third vertex gives a perfect code when $3\mid n$. Edges avoiding $C$ have cost zero. For an edge incident with $C$, orient and label the cycle so the deleted edge is $v_{3q}v_1$ and $v_1\in C$. The original code is $\{v_1,v_4,\ldots,v_{3q-2}\}$. By Lemma~\ref{l:pathcodes}, the resulting path has unique code $D=\{v_2,v_5,\ldots,v_{3q-1}\}$. The codes are disjoint and both have size $q$, so the cost is $2q$. Taking the maximum proves the result.
\end{proof}
\begin{corollary}\label{c:cycle}
	For every $q\geq1$, $R_{\min}(C_{3q})=R_{\max}(C_{3q})=2q$. In particular, there is no finite upper bound depending only on maximum degree, even for connected graphs with finite repair costs.
\end{corollary}
\begin{figure}[htbp]
	\centering
	\begin{tikzpicture}
		\foreach \s/\lab in {0/1,5/2}{
			\begin{scope}[xshift=\s cm]
				\foreach \i in {1,...,6}{\node[vertex] (v\lab\i) at ({90-60*(\i-1)}:1.2) {$v_{\i}$};}
				\foreach \i in {1,...,5}{\pgfmathtruncatemacro{\j}{\i+1}\draw[edge] (v\lab\i)--(v\lab\j);}
				\ifnum\lab=1
				\draw[failed] (v11)--(v16); \node[code] at (v11) {$v_1$}; \node[code] at (v14) {$v_4$};
				\node at (0,-1.9) {$C=\{v_1,v_4\}$};
				\else
				\node[code] at (v22) {$v_2$}; \node[code] at (v25) {$v_5$};
				\node at (0,-1.9) {$D=\{v_2,v_5\}$};
				\fi
		\end{scope}}
		\draw[-{Stealth},thick] (1.6,0)--(3.4,0) node[midway,above] {$-v_1v_6$};
	\end{tikzpicture}
	\caption{A cycle of order six in which all original code vertices must be replaced: $r(C_6,C;v_1v_6)=4$.}
\end{figure}
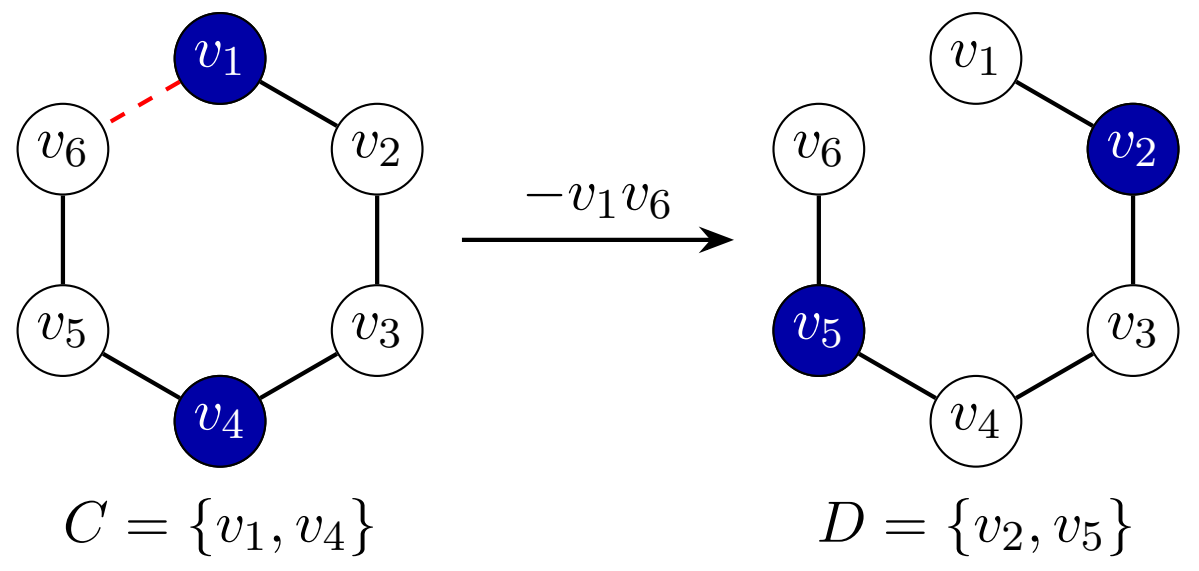

\subsection{Trees with irreparable failures}
\begin{example}\label{ex:tree}
	Let $T$ have vertices $a,b,c,d,e,f$ and edges $ab,bc,bd,ce,df$. The set $C=\{a,e,f\}$ is a perfect code. Deleting $ce$ isolates $e$, and the remaining component has leaves $a,c$ adjacent to $b$ and a path $bdf$. Every perfect code of that component must contain $b$: otherwise both $a$ and $c$ must belong to it, dominating $b$ twice. But selecting $b$ excludes $d$, and to dominate $f$ one must select $f$, which then dominates $d$ twice. Thus $T-ce$ has no perfect code. The automorphism interchanging $c$ with $d$ and $e$ with $f$ maps $ce$ to $df$, so $T-df$ also has no perfect code. These are the only two irreparable edges. Deleting $bc$ or $bd$ preserves $C$, and $\{a,c,f\}$ is a perfect code of $T-ab$. Hence $B(T)=\{ce,df\}$ and $R(T,C)=\infty$.
\end{example}
\begin{figure}[htbp]
	\centering
	\begin{tikzpicture}
		\foreach \s/\lab in {0/1,5/2}{
			\begin{scope}[xshift=\s cm]
				\node[vertex] (a\lab) at (0,1.6) {$a$}; \node[vertex] (b\lab) at (0,.5) {$b$};
				\node[vertex] (c\lab) at (-.9,-.5) {$c$}; \node[vertex] (d\lab) at (.9,-.5) {$d$};
				\node[vertex] (e\lab) at (-.9,-1.6) {$e$}; \node[vertex] (f\lab) at (.9,-1.6) {$f$};
				\draw[edge] (a\lab)--(b\lab)--(c\lab) (b\lab)--(d\lab)--(f\lab);
				\ifnum\lab=1 \draw[failed] (c1)--(e1); \node[code] at (a1) {$a$}; \node[code] at (e1) {$e$}; \node[code] at (f1) {$f$};
				\node at (0,-2.2) {$C=\{a,e,f\}$};
				\else \node at (0,-2.2) {No perfect code}; \fi
		\end{scope}}
		\draw[-{Stealth},thick] (1.4,0)--(3.6,0) node[midway,above] {$-ce$};
	\end{tikzpicture}
	\caption{A tree admitting a perfect code but having an irreparable edge failure.}
\end{figure}
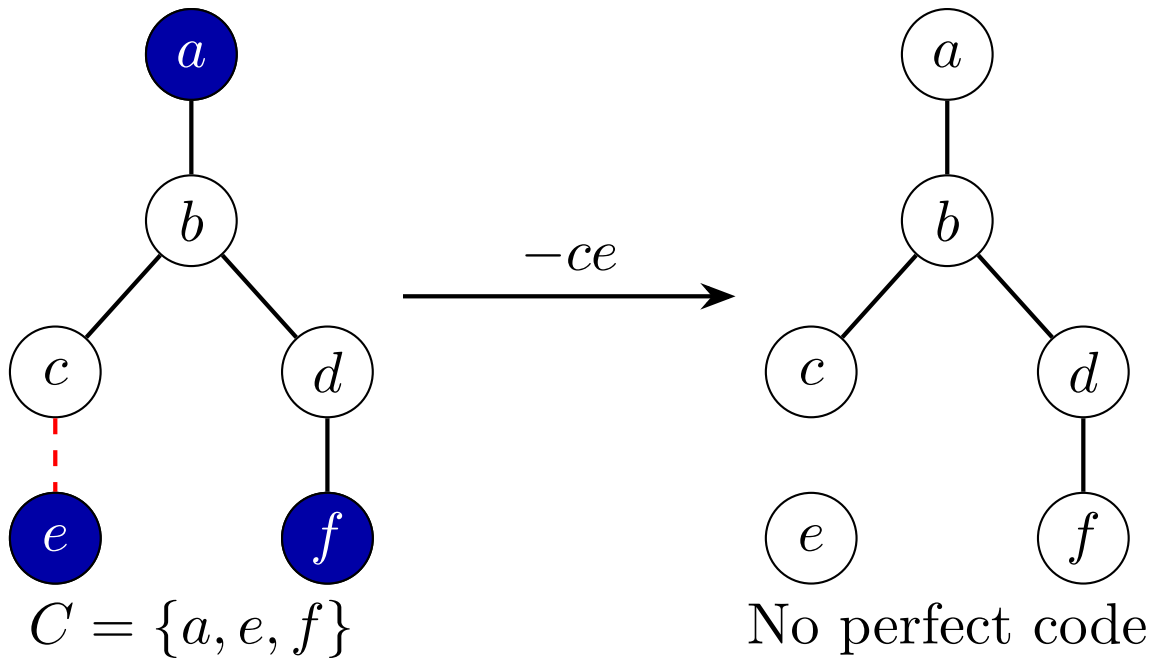

\section{Regular Graphs and Symmetry}
\subsection{Published results on domination and bondage}
We recall four results of Huang and Xu in the undirected setting. We state the parts used below without proof.

\begin{lemma}[{Huang and Xu~\cite[Lemma~2.3]{HuangXu2008}}]\label{hx:regular}
For a $d$-regular graph $\Gamma$ of order $n$, $\gamma(\Gamma)\ge n/(d+1)$. Equality holds precisely when $\Gamma$ admits a perfect code. In that case, the minimum dominating sets of $\Gamma$ are exactly its perfect codes.
\end{lemma}

An edge supports a dominating set $D$ if exactly one endpoint belongs to $D$. Let $s(\Gamma)$ be the least size of an edge set $F$ such that every minimum dominating set is supported by some edge in $F$.

\begin{lemma}[{Huang and Xu~\cite[Lemma~2.4]{HuangXu2008}}]\label{hx:support}
The bondage number satisfies $b(\Gamma)\ge s(\Gamma)$. If $\Gamma$ is regular and admits a perfect code, then $b(\Gamma)=s(\Gamma)$.
\end{lemma}

\begin{theorem}[{Huang and Xu~\cite[Theorem~3.1]{HuangXu2008}}]\label{hx:transitive}
If $\Gamma$ is an undirected vertex-transitive graph of positive degree and order $n$, then
\[
b(\Gamma)\ge \left\lceil\frac{n}{2\gamma(\Gamma)}\right\rceil .
\]
\end{theorem}

\begin{lemma}[{Huang and Xu~\cite[Lemma~4.3]{HuangXu2008}}]\label{hx:circulant}
Let $\Gamma$ be the circulant on $\mathbb Z_n$ with steps $\{\pm1,\pm s\}$, where $1<s<n/2$. Then $\Gamma$ admits a perfect code if and only if $5\mid n$ and $s\equiv\pm2\pmod 5$. Under these conditions its perfect codes are precisely the five residue classes modulo five.
\end{lemma}

\subsection{Consequences for edge repair}
Lemma~\ref{hx:regular} describes the minimum dominating sets before deletion. The next theorem identifies the perfect codes available after deletion.
\begin{theorem}\label{t:regular}
	Let $\Gamma$ be $d$-regular of order $n$ with $d\geq2$ and admitting a perfect code. For every edge $uv$,
	\[
	\PC(\Gamma-uv)=\{D\in\PC(\Gamma):u,v\notin D\}.
	\]
\end{theorem}
\begin{proof}
	Put $H=\Gamma-uv$. If $D$ is a perfect code of $H$, then counting the disjoint closed neighborhoods gives
	\[
	n=(d+1)|D|-|D\cap\{u,v\}|.
	\]
	Also $n=(d+1)|C|$. Subtracting gives $(d+1)(|D|-|C|)=|D\cap\{u,v\}|$. The right-hand side is an integer between zero and two, and $d+1\geq3$. It follows that both sides are zero, so $D$ avoids $u,v$. Adding the edge back cannot affect domination by $D$, giving $D\in\PC(\Gamma)$. Conversely, any perfect code of $\Gamma$ avoiding both endpoints survives the deletion.
\end{proof}
\begin{corollary}\label{c:regularbound}
Let $\Gamma$ be $d$-regular of order $n$, with $d\ge2$, and let $C\in\PC(\Gamma)$. For $e=uv$, put $\mathcal A_e=\{D\in\PC(\Gamma):D\cap\{u,v\}=\varnothing\}$. Then
\[
r(\Gamma,C;e)=
\begin{cases}
2\bigl(\gamma(\Gamma)-\max_{D\in\mathcal A_e}|C\cap D|\bigr),
   & \mathcal A_e\ne\varnothing,\\
\infty, & \mathcal A_e=\varnothing.
\end{cases}
\]
Every finite repair cost is even. If $R(\Gamma,C)<\infty$, then
\[
2\le R(\Gamma,C)\le 2\gamma(\Gamma)=\frac{2n}{d+1}.
\]
\end{corollary}
\begin{proof}
By Theorem~\ref{t:regular}, $\PC(\Gamma-e)=\mathcal A_e$. If this family is empty, the cost is infinite by definition. Otherwise each $D\in\mathcal A_e$ has the same cardinality as $C$, namely $\gamma(\Gamma)=n/(d+1)$. Thus $|C\dist D|=2(\gamma(\Gamma)-|C\cap D|)$, and minimizing over $D$ gives the formula. The upper bound and parity follow immediately. Since $d\ge2$, some edge meets $C$ and has positive cost, so a finite maximum is at least two. Cycles of order divisible by three and complete graphs of order at least three attain the upper bound.
\end{proof}
\begin{corollary}\label{c:regular-bondage}
	For a $d$-regular graph with $d\geq2$ admitting a perfect code,
	\[
	B(\Gamma)=\bigcap_{C\in\PC(\Gamma)} E_C(\Gamma).
	\]
	In particular, $R_{\max}(\Gamma)<\infty$ if and only if this intersection is empty.
\end{corollary}
\begin{proof}
	By Theorem~\ref{t:regular}, an edge is irreparable precisely when no perfect code of $\Gamma$ avoids both endpoints.
\end{proof}
Suppose $\Gamma$ is $d$-regular with $d\ge2$ and admits a perfect code. By Lemmas~\ref{hx:regular} and~\ref{hx:support}, $b(\Gamma)=1$ exactly when a single edge supports every perfect code. Since perfect codes are independent, Theorem~\ref{t:regular} yields
\[
 b(\Gamma)=1 \quad\Longleftrightarrow\quad B(\Gamma)\ne\varnothing .
\]
Indeed, a single edge supports every perfect code precisely when no perfect code avoids both endpoints.

For a vertex-transitive graph of degree $d$ admitting a perfect code, Theorem~\ref{hx:transitive} and Lemma~\ref{hx:regular} give $b(\Gamma)\ge\lceil(d+1)/2\rceil$, the lower bound stated explicitly in Huang and Xu~\cite[Corollary~3.4]{HuangXu2008}. For $d\ge2$ the bound excludes a single irreparable edge by Corollary~\ref{c:regular-bondage}. This yields the degree-at-least-two case of the next theorem from their results. We include a direct proof that also identifies a repair in the automorphism orbit of the original code.
\begin{theorem}\label{t:vt}
	Every finite vertex-transitive graph of positive degree admitting a perfect code has finite maximum edge repair costs.
\end{theorem}
\begin{proof}
	A vertex-transitive graph is regular, say of degree $d$. If $d=1$, it is a disjoint union of copies of $K_2$, and every edge has repair cost one. Suppose $d\geq2$. Fix a perfect code $C$, choose an automorphism $\varphi$ uniformly at random, and let $uv$ be any edge. By transitivity, each vertex belongs to $\varphi(C)$ with probability $|C|/n=1/(d+1)$. Since perfect codes are independent, the events that $u$ or $v$ belongs to $\varphi(C)$ are disjoint. Consequently,
	\[
	\Pr\bigl(u,v\notin\varphi(C)\bigr)=1-\frac{2}{d+1}>0.
	\]
	There is therefore a perfect code avoiding both endpoints of the edge. It remains a perfect code after the deletion. This holds for every edge, so Proposition~\ref{p:finite} yields the conclusion.
\end{proof}
\begin{corollary}\label{c:cayley}
	Every finite undirected Cayley graph of positive degree admitting a perfect code has finite maximum edge repair costs.
\end{corollary}
\begin{proof}
	Every Cayley graph is vertex-transitive; apply Theorem~\ref{t:vt}.
\end{proof}

Lemma~\ref{hx:circulant} and Theorem~\ref{t:regular} also give an exact repair cost for double-loop circulants.

\begin{corollary}\label{c:circulant}
Let $\Gamma$ have vertex set $\mathbb Z_n$ and steps $\{\pm1,\pm s\}$, where $1<s<n/2$. If $5\mid n$ and $s\equiv\pm2\pmod5$, then every perfect code $C$ satisfies $R(\Gamma,C)=2n/5$.
\end{corollary}
\begin{proof}
By Lemma~\ref{hx:circulant}, the perfect codes are the five residue classes modulo five. For an edge incident with $C$, Theorem~\ref{t:regular} permits exactly those classes avoiding both endpoints. At least three such classes remain, each disjoint from $C$ and of size $n/5$. The repair cost is therefore $2n/5$. All edges avoiding $C$ have cost zero.
\end{proof}
\section{Disjoint Unions, Automorphisms, and Extremal Questions}
\begin{proposition}\label{p:union}
	Let $\Gamma_1,\Gamma_2$ be vertex-disjoint graphs, each with at least one edge, and $C_i\in\PC(\Gamma_i)$. If $\Gamma=\Gamma_1\cup \Gamma_2$, then
	\[
	R(\Gamma,C_1\cup C_2)=\max\{R(\Gamma_1,C_1),R(\Gamma_2,C_2)\}.
	\]
\end{proposition}
\begin{proof}
	For an edge $e$ of $\Gamma_1$, any perfect code after deletion decomposes as $D_1\cup D_2$, where $D_1\in\PC(\Gamma_1-e)$ and $D_2\in\PC(\Gamma_2)$. The symmetric-difference cost is the sum of the costs on the two components. Choosing $D_2=C_2$ minimizes its contribution, so the minimum equals $r(\Gamma_1,C_1;e)$, including the case when $\PC(\Gamma_1-e)$ is empty. The same argument applies to edges of $\Gamma_2$. Taking the maximum proves the assertion.
\end{proof}
\begin{corollary}\label{c:union}
	Under the hypotheses of Proposition~\ref{p:union},
	\begin{align*}
		R_{\min}(\Gamma_1\cup \Gamma_2)&=\max\{R_{\min}(\Gamma_1),R_{\min}(\Gamma_2)\},\\
		R_{\max}(\Gamma_1\cup \Gamma_2)&=\max\{R_{\max}(\Gamma_1),R_{\max}(\Gamma_2)\}.
	\end{align*}
\end{corollary}
\begin{proof}
	Every perfect code of the union is the union of component codes. Minimizers can be chosen independently on the two components, and likewise for maximizers. Proposition~\ref{p:union} gives both identities.
\end{proof}
\begin{remark}\label{rem:isolates}
	Isolated vertices can be omitted because they belong to every perfect code and never change their membership following an edge deletion.
\end{remark}
\begin{proposition}\label{p:automorphism}
	For any automorphism $\varphi$ of $\Gamma$ and any perfect code $C$,
	\[
	R(\Gamma,C)=R(\Gamma,\varphi(C)).
	\]
\end{proposition}
\begin{proof}
	The map $\varphi$ sends the perfect codes of $\Gamma-e$ bijectively to the perfect codes of $\Gamma-\varphi(e)$, and preserves symmetric-difference cardinalities. Thus $r(\Gamma,C;e)=r(\Gamma,\varphi(C);\varphi(e))$, including infinite costs. Maximizing over the edges gives the result.
\end{proof}
\begin{corollary}\label{c:code-transitive}
	If the automorphism group of $\Gamma$ acts transitively on its family of perfect codes, then $R_{\min}(\Gamma)=R_{\max}(\Gamma)$.
\end{corollary}

\begin{definition}\label{def:extremal}
	For $n\geq2$, define
	\[
	M(n)=\max\{R(\Gamma,C): |V(\Gamma)|=n,\ \Gamma\text{ connected},\ C\in\PC(\Gamma),\ R(\Gamma,C)<\infty\}.
	\]
\end{definition}
\begin{proposition}\label{p:extremal}
	For every $n\geq2$,
	\[
	\left\lceil\frac{2n-1}{3}\right\rceil\leq M(n)\leq n-1.
	\]
	In particular, $M(3q+1)\ge 2q+1$ for every $q\ge1$.
\end{proposition}
\begin{proof}
	For $n=3q$, the cycle gives cost $2q$ by Theorem~\ref{t:cycle}. For $n=3q+1$, use Theorem~\ref{t:pendant-cycle}. For $n=3q+2$, including $n=2$, the path gives cost $2q+1$ by Theorem~\ref{t:path}. Proposition~\ref{p:bounds} gives the upper bound.
\end{proof}
The exhaustive census in Appendix~\ref{app:census} gives the exact connected maxima
\[
\begin{array}{c|rrrrrr}
n&2&3&4&5&6&7\\\hline
M(n)&1&2&3&3&4&5
\end{array}
\]
Thus $M(n)=n-2$ for $5\le n\le7$. The program includes graphs with a unique perfect code and maximizes only over connected graphs with finite costs.

\begin{problem}\label{prob:extremal}
	Determine $M(n)$ for $n\ge8$ and characterize the pairs $(\Gamma,C)$ attaining it.
\end{problem}
The dependence on the initial code in Example~\ref{ex:depend} persists for arbitrarily long tails. Count tail length in edges, starting at the attachment vertex.
\begin{proposition}\label{p:clique-tail-spread}
Let $m\ge4$ and $k\ge1$. Attach the path $x t_1\cdots t_{3k}$ to $K_m$ at $x$, and denote the resulting graph by $\Gamma$. Then
\[
R_{\min}(\Gamma)=2k,\qquad R_{\max}(\Gamma)=2k+2,\qquad
\sigma_R(\Gamma)=2.
\]
\end{proposition}
\begin{proof}
Exactly one clique vertex belongs to a perfect code: two would be adjacent, and without one the clique vertices other than $x$ would be uncovered. Fix $b\in V(K_m)\setminus\{x\}$. On the path $bxt_1\cdots t_{3k}$, every perfect code contains exactly one of $b,x$, since these are the only vertices covering $b$. Such a code extends unchanged to the clique, as each omitted clique vertex is adjacent to both $b$ and $x$. Conversely, every code whose clique vertex is $b$ or $x$ restricts to a code of this path. By Lemma~\ref{l:pathcodes}, the codes of $\Gamma$ are therefore
\[
 C_x=\{x\}\cup\{t_{3j}:1\le j\le k\},\qquad
 C_b=\{b\}\cup\{t_{3j-1}:1\le j\le k\}
 \quad(b\in V(K_m)\setminus\{x\}).
\]

First consider $C_b$. A clique-edge deletion costs at most two: choose a nonattachment clique vertex outside the two endpoints as the clique center, retaining all tail centers. Tail edges disjoint from $C_b$ cost zero. For an edge incident with $t_{3j-1}$, retain the original centers on the component containing this center. If the cut is just after it, the uncovered side is a path of order $3(k-j)+1$. Its unique code has $k-j+1$ vertices, replacing at most $k-j$ old centers, for cost at most $2(k-j)+1$. If the cut is just before it, the other component consists of the clique and $t_1,\ldots,t_{3j-2}$. Its unique code is $\{x,t_3,t_6,\ldots,t_{3j-3}\}$, where the tail part is empty for $j=1$. Indeed, restricting to the path starting at any candidate nonattachment clique center as above gives a path of order $3j$, whose unique code starts at $x$. The $j$ new centers are disjoint from the $j$ old ones on this component, so the cost is exactly $2j$. Taking $j=k$ gives $R(\Gamma,C_b)=2k$.

Now consider $C_x$. After a tail-edge deletion, any repaired code on $bxt_1\cdots t_{3k}$ still contains exactly one of $b,x$ and hence extends unchanged to the full graph. Theorem~\ref{t:path} on this path of order $3k+2$ gives cost at most $2k+1$. After any clique-edge deletion some $C_b$ avoids both endpoints, giving cost at most $2k+2$. To see equality, delete $xb$. Every repaired code contains exactly one clique vertex: a third clique vertex requires coverage, and even the nonadjacent pair $x,b$ would cover that third vertex twice. Selecting $x$ would leave $b$ uncovered. Selecting $b$ would require $t_1$ to cover $x$, yet the tail would then have to be covered internally, and the unique code of $P_{3k}$ omits $t_1$. Thus the clique center is some $b'\notin\{x,b\}$ and the tail has its unique internal code. All such repairs are disjoint from $C_x$ and have size $k+1$, proving $R(\Gamma,C_x)=2k+2$.
\end{proof}

For $m=4$, the first three pairs of extrema are $(2,4)$, $(4,6)$, and $(6,8)$. The two extrema grow with the tail length while their difference stays constant. Taking $k=1$ and varying $m$ also realizes spread two at every order $n\ge7$.

\begin{problem}\label{prob:spread}
	Characterize the connected graphs with $R_{\min}(\Gamma)=R_{\max}(\Gamma)<\infty$, and determine whether $\sigma_R(\Gamma)$ is unbounded among connected fully edge-repairable graphs.
\end{problem}
Every positive integer occurs as the maximum repair cost of a perfect code in a tree. Indeed, for $k\ge1$ put $n=\lceil(3k+1)/2\rceil$. Then $3k\le 2n-1\le 3k+1$, so Theorem~\ref{t:path} gives $R(P_n,C)=k$ for every perfect code $C$ of $P_n$.

\begin{problem}\label{prob:trees}
	Characterize the trees $T$ with at least one edge that admit a perfect code and for which $T-e$ admits a perfect code for every $e\in E(T)$.
\end{problem}
\begin{problem}\label{prob:regular}
	Determine which regular graphs admitting perfect codes fail to be fully edge-repairable, and characterize the extremal cases in Corollary~\ref{c:regularbound}.
\end{problem}
\begin{table}[htbp]
	\centering\small
	\caption{Selected consequences of the additional results.}
	\begin{tabular}{p{.46\linewidth}p{.43\linewidth}}
		\toprule
		Condition or family & Result\\\midrule
		Finite maximum repair cost & $1\leq R(\Gamma,C)\leq\min\{n-1,2\gamma(\Gamma)+1\}$\\
		Repairable edge & $r(\Gamma,C;e)\equiv\gamma(\Gamma-e)-\gamma(\Gamma)\pmod2$\\
		Star forest with an edge & $R_{\min}=R_{\max}=1$\\
		Complete graph, $n\geq3$ & $R_{\min}=R_{\max}=2$\\
		Path, $n\geq2$ & $R_{\min}=R_{\max}=\lfloor(2n-1)/3\rfloor$\\
		Cycle, $n=3q$ & $R_{\min}=R_{\max}=2q$\\
		Regular graph, $d\geq2$, finite costs & $2\leq R(\Gamma,C)\leq 2n/(d+1)$\\
		Vertex-transitive, positive degree & $R_{\max}<\infty$ whenever a code exists\\
		\bottomrule
	\end{tabular}
\end{table}
\noindent The domination-number identity and parity proposition connect repair costs to classical domination theory. The regular-graph characterization identifies all repairs, and the vertex-transitive argument locates a surviving code in the initial orbit when the degree is at least two. The path and cycle formulas show that finite repair costs can grow linearly with graph order even at bounded degree.

\subsection{Repair models and computational problems}
Let $C\in\PC(\Gamma)$. In each repair model below, the minimum over an empty family is $\infty$.

\paragraph{Simultaneous edge failures.}
For $1\le k\le |E(\Gamma)|$ and $F\subseteq E(\Gamma)$ with $|F|=k$, define the repair cost and its worst-case value by
\[
r_k(\Gamma,C;F)=\min_{D\in\PC(\Gamma-F)}|C\dist D|,
\qquad
R_k(\Gamma,C)=\max_{\substack{F\subseteq E(\Gamma)\\|F|=k}}r_k(\Gamma,C;F).
\]
The graph is fully repairable under $k$ simultaneous edge failures if $\PC(\Gamma-F)\ne\varnothing$ for every such $F$.

\paragraph{Single-edge insertion.}
For distinct nonadjacent vertices $u,v$, the insertion repair cost is
\[
r^{+}(\Gamma,C;uv)=\min_{D\in\PC(\Gamma+uv)}|C\dist D|.
\]

\paragraph{Single-vertex failure.}
For $x\in V(\Gamma)$, define the repair cost on the surviving vertices by
\[
r_{\mathrm v}(\Gamma,C;x)=
\min_{D\in\PC(\Gamma-x)}|(C\setminus\{x\})\dist D|.
\]
Under this convention the forced loss of $x$ incurs no membership-change charge.

\paragraph{Computational problems.}
The single-edge repair-cost function problem takes as input a graph $\Gamma$, a specified perfect code $C$, and an edge $e$, and returns $r(\Gamma,C;e)$, including the value $\infty$. Its budget decision version takes an additional nonnegative integer $h$ and asks whether $r(\Gamma,C;e)\le h$. The full edge-repairability decision problem, on input a graph $\Gamma$ promised to have at least one edge and a perfect code, asks whether $\PC(\Gamma-e)\ne\varnothing$ for every $e\in E(\Gamma)$.

\paragraph{Dominating-set reconfiguration.}
For a fixed graph $H$ and an integer $b\ge0$, the token addition and removal reconfiguration graph has as vertices the dominating sets of $H$ of size at most $b$. Two sets are adjacent when their symmetric difference has size one. The reconfiguration distance is the length of a shortest path between two such sets, or $\infty$ if no path exists. In the repair model, the distance charged between $C\in\PC(\Gamma)$ and $D\in\PC(\Gamma-e)$ is the endpoint distance $|C\dist D|$. Intermediate sets are unrestricted, and $C$ may already fail to dominate $\Gamma-e$.

\subsection*{Conclusion}
The direct proof of Theorem~\ref{t:vt} describes how symmetry supplies a repair. If a finite vertex-transitive graph $\Gamma$ has positive degree and admits a perfect code, then $\Gamma-e$ admits a perfect code for every $e\in E(\Gamma)$. For degree at least two, an automorphic image of the original code can be chosen to avoid both endpoints of the failed edge and therefore survives its deletion. In degree one, deleting an edge isolates its endpoints, and adding the previously unselected endpoint repairs the code. The cost of repair can nevertheless grow with the order of the graph: on $C_{3q}$ the worst-case cost is $2q$ by Theorem~\ref{t:cycle}.

We have established general constraints on single-edge repair costs and derived exact values for several fundamental graph families. The results show that the cost can depend on the initial perfect code and can grow linearly even in graphs of bounded degree. Further work should identify additional families with finite repair costs, sharpen the general bounds, and complete the comparison with the literature on efficient domination and edge-deletion vulnerability.
\appendix
\section{Reproducing the small-order census}\label{app:census}
Save the listing as \texttt{census.js} and run \texttt{node census.js} with Node.js. No additional packages are required. Each output row reports the order, the connected maximum $M(n)$, and a dictionary mapping $(m,R_{\min},R_{\max})$ to its labelled count among graphs with unequal extrema.

The array $N$ stores closed neighborhoods as bit masks. A perfect code partitions the vertex set into closed neighborhoods of its members. At each recursive step, the first uncovered vertex $v$ must belong to exactly one block of the final partition. The procedure tries each possible center $u\in N[v]$ whose block is disjoint from the blocks already chosen. Hence each perfect code reaches exactly one completed branch. The same procedure enumerates repairs after deleting each edge. Exclusive OR computes symmetric difference, and the minimum of an empty list is positive infinity in JavaScript.
{\small
\begin{verbatim}
// Run with Node.js: node census.js
function codes(N) {
  const all = (1 << N.length) - 1, found = [];
  function visit(covered, chosen) {
    if (covered === all) { found.push(chosen); return; }
    const missing = all ^ covered;
    const v = 31 - Math.clz32(missing & -missing);
    for (let u = 0; u < N.length; ++u) {
      if ((N[v] & (1 << u)) && !(covered & N[u]))
        visit(covered | N[u], chosen | (1 << u));
    }
  }
  visit(0, 0);
  return found;
}
function census(n) {
  const edges = [], counts = new Map();
  let maximum = 0;
  const size = Array(1 << n).fill(0);
  for (let i = 1; i < size.length; ++i)
    size[i] = 1 + size[i & (i - 1)];
  for (let u = 0; u < n; ++u)
    for (let v = u + 1; v < n; ++v) edges.push([u, v]);
  for (let mask = 1; mask < (1 << edges.length); ++mask) {
    const N = Array.from({length:n}, (_, u) => 1 << u);
    const present = edges.filter((_, i) => mask & (1 << i));
    for (const [u,v] of present) {
      N[u] |= 1 << v; N[v] |= 1 << u;
    }
    const P = codes(N);
    if (!P.length) continue;
    const worst = P.map(() => 0);
    for (const [u,v] of present) {
      N[u] ^= 1 << v; N[v] ^= 1 << u;
      const Q = codes(N);
      N[u] ^= 1 << v; N[v] ^= 1 << u;
      for (let i = 0; i < P.length; ++i) {
        const cost = Math.min(...Q.map(d => size[P[i] ^ d]));
        worst[i] = Math.max(worst[i], cost);
      }
      if (!Q.length) break;
    }
    const lo = Math.min(...worst), hi = Math.max(...worst);
    if (Number.isFinite(hi) && hi > maximum) {
      let reached = 1, old;
      do {
        old = reached;
        for (let u = 0; u < n; ++u)
          if (reached & (1 << u)) reached |= N[u];
      } while (reached !== old);
      if (reached === (1 << n) - 1) maximum = hi;
    }
    if (lo < hi) {
      const key = [present.length, lo, hi].join(",");
      counts.set(key, (counts.get(key) || 0) + 1);
    }
  }
  console.log(n, maximum, Object.fromEntries(counts));
}
for (let n = 2; n <= 7; ++n) census(n);
\end{verbatim}
}

\end{document}